%% file: agt-2-13.tex
\documentclass{gtart}

\input agtout

\lognumber{13}
\volumenumber{2}
\volumeyear{2002}
\papernumber{13}
\published{18 April 2002}
\pagenumbers{285}{296}
\received{6 December 2001}
\accepted{10 April 2002}

\newtheorem{theorem}{Theorem}[section]    
 \newtheorem{lemma}[theorem]{Lemma}        
 
 \newtheorem*{theoremrepeat}{Theorem \ref{introflatinject}}

 \def\SL{\mbox{\rm{SL}}} \def\SO{\mbox{\rm{SO}}} \def\O{\mbox{\rm{O}}}
 \def\Isom{\mbox{Isom}} 
  \def\GL{\mbox{\rm{GL}}} 
 \def\O{\mbox{\rm{O}}} \def\A{\cal A}
  \def\stab{\mbox{\rm{Stab}}}

 \def\demo{\proof}

\begin{document}
\title{All flat manifolds are cusps of hyperbolic orbifolds} 
\authors{D.D. Long\\A.W. Reid}                  
\address{Department of Mathematics, University of California\\Santa Barbara,
CA 93106, USA}                  
\secondaddress{Department of Mathematics, University of Texas\\Austin, 
TX 78712, USA}             
\asciiaddress{Department of Mathematics, University of
California\\Santa Barbara, CA 93106, USA\\and\\Department of
Mathematics, University of Texas\\Austin, TX 78712, USA}

\email{long@math.ucsb.edu, areid@math.utexas.edu}

\begin{abstract} 

We show that all  closed flat $n$-manifolds are diffeomorphic
to a cusp cross-section  in a  finite volume hyperbolic $n+1$-orbifold.

\end{abstract}
\asciiabstract{We show that all closed flat n-manifolds are
diffeomorphic to a cusp cross-section in a finite volume hyperbolic
(n+1)-orbifold.}


\maketitle 
\let\\\par

 \section{Introduction}

 By a flat {\em n-manifold} (resp.\ flat {\em n-orbifold}) we mean a
 manifold (resp.\ orbifold) ${\bf E}^n/\Gamma$ where ${\bf E}^n$ is
 Euclidean n-space and $\Gamma$ a discrete, cocompact torsion-free
 subgroup of $\Isom ({\bf E}^n)$ (resp.\ $\Gamma$ has elements of finite
 order). In \cite{HR},  Hamrick and Royster resolved a longstanding
 conjecture by showing that {\em every} flat
 $n$-manifold bounds an $(n+1)$-dimensional manifold, in the sense that
 each diffeomorphism class has a representative that bounds.  

 Flat manifolds and orbifolds are connected with hyperbolic
 orbifolds via the structure of the
 cusp ends of finite volume hyperbolic orbifolds---if $M^{n+1}$ is a
 non-compact finite volume hyperbolic orbifold, then by a standard
 analysis of the thin parts, a cusp cross-section is a flat
 $n$-orbifold (see below and \cite{Th}).  In the early eighties, motivated
 by this and work of Gromov \cite{G}, Farrell and Zdravkovska \cite{FZ}
 posed the following geometric version:\\
 \vspace{2pt}
 {\bf Question}\qua Does the
 diffeomorphism class of every flat $n$-manifold have a representative
 $W$ which arises as the cusp cross-section of a finite volume
 $1$-cusped hyperbolic $(n+1)$-manifold?\\
 \vspace{2pt}
 We note that
 it makes sense to ask only for bounding up to diffeomorphism type,
 since it is a well-known consequence of Mostow Rigidity that there are
 algebraic restrictions on the {\em isometry type} of the flat
 $n$-manifolds that can arise as cusp cross sections of finite volume
 hyperbolic $(n+1)$-manifolds.

 In this generality, the question has a positive answer in dimension $2$,
 that is, for the torus and Klein bottle, but as is shown
 in \cite{LR} is false in dimension 3 (and indeed in all
dimensions $4k-1 \geq 3$).
 On the other hand, since finite volume, non-compact hyperbolic
 manifolds exist in all dimensions (see \S 6 for example), this
 together with the fact that maximal abelian subgroups are separable
 (see below and \cite{Lo}) can be used to deduce the $n$-torus is a
 cusp cross section of a hyperbolic manifold with possibly many cusps.

 In this paper we shall show:
\begin{theorem}
\label{introgromov} For every $ n \geq 2$, the diffeomorphism class of
every flat $n$-manifold has a representative $W$ which arises as some
cusp cross-section of a finite volume cusped hyperbolic
$(n+1)$-orbifold.
 \end{theorem}
 
 As already remarked, much more is known in dimension $2$, and
 we only include this since the argument also works there.
 In dimension $3$, \cite{Ni}  proves a little more, namely that every
 flat 3-manifold is a  cusp cross-section of a hyperbolic 4-{\em manifold}, but
 even there, the number of cusps is not known to be one. Also
 in dimension 3, in \cite{RT}, it is shown by an {\em ad hoc} argument
 that of the $10$ diffeomorphism types of flat 3-manifolds,
 $7$ have representatives that arise as a cusp cross-section in some 5 cusped
 hyperbolic manifold related to gluings of the ideal $24$-cell in ${\bf H}^4$.

 The proof of Theorem \ref{introgromov} uses arithmetic methods, together with
 some extra control one can arrange to make a separability argument.  The main
 ingredient
 in the proof is the following (see \S 2 for definitions and notation):

 \begin{theorem} \label{introflatinject}
Let $\Gamma$ be the fundamental group  of a flat $n$-manifold.
 Then there is a quadratic form $q_{n+2}$ defined over $\bf Q$, of
 signature $(n+1,1)$ for which $\Gamma$ embeds as a subgroup of
 $\O_0(q_{n+2};{\bf Z})$. \end{theorem}

The group $\O_0(q_{n+2};{\bf Z})$ is arithmetic, and hence of
 finite co-volume acting on ${\bf H}^{n+1}$.
 The proof of Theorem \ref{introgromov} is then completed by using a
 subgroup separability argument to pass to a subgroup of finite index
 in $\O_0(q_{n+2};{\bf Z})$ for which the group $\Gamma$ is a maximal
 peripheral subgroup.

 As with \cite{Ni}, we cannot guarantee 1-cusped examples,
 and at present, we have been unable to pass to manifolds, even  
 assuming separability of geometrically finite subgroups.

 It appears that the question in its full generality is much harder to
 approach, for example, it appears be unknown, whether,  for
 $n \geq 4$, there even {\em exist} 1-cusped hyperbolic $n$-manifolds of finite
 volume. As far as the authors are aware, this is known only in dimensions $2$
 and $3$.

 \section{Preliminaries.}
 \label{preliminaries}

 \subsection{Quadratic forms}

 We need to recall some standard facts about quadratic forms and
 orthogonal groups of such forms; \cite{La} is a standard reference.

 If $f$ is a quadratic form in $n+1$ variables with
 coefficients in $K$ and associated symmetric matrix $F$, let $$\O(f) = \{X
 \in \GL(n+1,{\bf C})~|~X^tFX = F \}$$ be the {\em Orthogonal group} of $f$,
 and $$\SO(f) = \O(f) \cap \SL(n+1,{\bf C}),$$ the {\em Special Orthogonal
 group} of $f$.  These are algebraic groups defined over $K$. If $L$ is
 a subring of $\bf C$ we denote the set of $L$-points of these groups
 by $\O(f;L)$ (resp.\ $\SO(f;L)$).\\
 \vspace{2pt}
 {\bf Definition.} Two $n$-dimensional quadratic forms
 $f$ and $q$ defined over a field $K$ (with associated symmetric matrices
 $F$ and $Q$) are {\em equivalent} over $K$ if there exists $P \in \GL(n,
 K)$ with $P^tFP = Q$.\\
 \vspace{2pt}
 If $K\subset {\bf R}$ is a number
 field, and $R_K$ its ring of integers, then $\SO(f;R_K)$ is an arithmetic
 subgroup of $\SO(f;{\bf R})$, \cite{BH} or \cite{Bo}. In particular
 $\SO(f;R_K)$ has finite co-volume acting on an associated
 symmetric space. The following is well-known and proved in \cite{ALR} for
 example.

 \begin{lemma} \label{equivalent} If $f$ and $q$ are equivalent over
$K$ then:
\begin{itemize}
 \item $\SO(f;{\bf R})$ is conjugate to $\SO(q;{\bf R})$ and $\SO(f;K)$ is
 conjugate to\nl $\SO(q;K)$. \item $\SO(f;R_K)$ is conjugate to a subgroup of
 $\SO(q;K)$ commensurable with $\SO(q;R_K)$.\qed \end{itemize}
 \end{lemma}

 \subsection{Crystallographic groups and hyperbolic orbifolds}

 We record some facts about crystallographic and Bieberbach groups that
 we will need, see
 \cite{Ch} for a comprehensive discussion of these groups, and \cite{Th}
 Theorem 4.2.2 for the theorem stated below.

 An {\em n-dimensional crystallographic group } is a cocompact
 discrete group of isometries of ${\bf E}^n$.
 When $\Gamma$ is torsion-free it is called a {\em Bieberbach group}.
 By Bieberbach's Third Theorem (see \cite{Ch}),
 the number of n-dimensional
 crystallographic groups up to affine equivalence is finite.
 What we require is summarized in the following from \cite{Th} page 222:

\begin{theorem}
 \label{crystal}
An $n$-dimensional crystallographic group $\Gamma$ contains a normal
subgroup of finite index consisting of translations, that is free
abelian of rank $n$.  The maximal such subgroup is characterized as the unique
maximal abelian subgroup of finite index in $\Gamma$.\qed\end{theorem}

 Theorem \ref{crystal} implies that associated to an $n$-dimensional
 crystallographic group $\Gamma$ is a finite group $\theta (\Gamma)$, the
 holonomy group, and an extension:
 $$1 \longrightarrow {\bf Z}^n \longrightarrow \Gamma
                       \longrightarrow \theta (\Gamma)\longrightarrow 1.$$
 When $\Gamma$ is a Bieberbach group we get a free action on ${\bf E}^n$,
 by rigid motions, that is for all $\gamma \in \Gamma$, and $v \in {\bf E}^n$,
 $$v \rightarrow \theta(\gamma )v + t_\gamma,$$
 for some $t_\gamma \in {\bf E}^n$ (see \cite{Ch} for details).

 \subsubsection{}

 We refer the reader to  \cite{Ra} or \cite{Th} for further details
 on what is contained in the next two subsections.

 Equip ${\bf R}^{n+1}$ with the quadratic form $f_n = \langle -1,1\ldots
 ,1,1\rangle$ of signature $(n,1)$.  The connected component of the identity
 in $\O(f_n;\!{\bf R})$ will be denoted $\O_0(f_n;\!{\bf R})$. This group
 preserves the upper sheet of the hyperboloid $f_n(x) = -1$ but
 contains reflections so reverses orientation.  We identify
 $\O_0(f_n;{\bf R})$ with $\Isom ({\bf H}^n)$. Passing to the
 connected component of the identity in $\SO(f_n;{\bf R})$, denoted
 $\SO_0(f_n;{\bf R})$ (which has index $4$ in $\O(f_n;{\bf R})$),
 gives a group which may be identified with $\Isom_+({\bf H}^n)$; it
 preserves the upper sheet of the hyperboloid $f_n(x) = -1$ and the
 orientation.  Given a (discrete) subgroup $\Delta$ of $\O(n,1;{\bf
 R})$, $\Delta \cap \SO_0(n,1;{\bf R})$ has index $\leq 4$ in
 $\Gamma$.

 \subsubsection{} An element of $\O_0(f_n;{\bf R})$ is parabolic
(resp.\ elliptic) if it has
 a unique fixed point which lies on $S^{n-1}_\infty$ (resp.\ has a fixed point
in ${\bf H}^n$).
Given a non-compact hyperbolic
 $n$-orbifold $Q = {\bf H}^n/\Delta$ of finite volume, and $C$ a cusp
cross-section
of $Q$, then there is a subgroup $\Delta_C <  \Delta$
consisting of parabolic and elliptic
elements having an
invariant horosphere $\cal H$ such that ${\cal H}/\Delta_C = C$. The group
$\Delta_C$ is a crystallographic group. This group is called a maximal
peripheral subgroup of $\Delta$.

In terms of the model above, an element is parabolic if and only
if it is not elliptic and leaves invariant a unique
light-like vector $v$. Furthermore, in the context of $\Delta$ and $\Delta_C$
above, all elements of $\Delta_C$ will preserve this unique light-like
vector. We summarize what we need.

 \begin{lemma}
 \label{peripheral}
 Let $Q={\bf H}^{n+1}/\Delta$ be a non-compact finite volume hyperbolic
 orbifold. $\Delta_C$ is a maximal peripheral subgroup of
 $\Delta$ if and only if $\Delta_C$ is the maximal subgroup of $\Delta$
 leaving the light-like vector $v$ invariant. Furthermore when
 $\Delta_C$ is maximal, by choice
 of a horosphere $\cal H$, ${\cal H}/\Delta_C  \cong C$ is an embedded cusp
 cross-section of $Q$.\qed\end{lemma}

 \subsubsection{}
 We record the following
 for convenience concerning arithmetic subgroups of
 $\Isom({\bf H}^n)$. For more details, see \cite{BH}, \cite{Bo} and
 \cite{VS}.

 Let $f$ be a diagonal quadratic form with rational coefficients and signature
 $(n,1)$. Thus there is a $P\in \GL(n+1,{\bf R})$ such that
 $P^tFP = F_n$, and so
 the group $P\O_0(f;{\bf Z})P^{-1}$ defines a discrete
 arithmetic subgroup of $\Isom ({\bf H}^n)$. The theory of arithmetic
 groups then gives,

 \begin{theorem} \label{noncocompact}
 In the notation above, the groups $P\O_0(f;{\bf Z})P^{-1}$ are finite
 co-volume subgroups of $\Isom ({\bf H}^n)$.
 \qed\end{theorem}

 In what follows, we will suppress the conjugation by $P$.
 A group $\O_0(f;{\bf Z})$ (and hence the conjugate in $\Isom({\bf
 H}^n)$) is cocompact if and only if the form $f$ does not represent $0$
 non-trivially with values in $\bf Q$, see \cite{BH}.  Whenever $n \geq 4$, the
 arithmetic groups constructed above are non-cocompact,
 since it is well known every indefinite quadratic
 form over $\bf Q$ in at least $5$ variables represents $0$ non-trivially,
 see \cite{La}.
 In fact up to commensurability, all non-cocompact
 arithmetic subgroups of $\O_0(f_n;{\bf R})$ arise from this construction
 (see \cite{VS}).

 \subsection{Some technical lemmas}

 In this section we gather together a collection of well-known
 results on separability properties of groups that will be used to pass from
 Theorem \ref{introflatinject} to Theorem \ref{introgromov}.

 Recall that a subgroup $H$ of a group $G$ is {\em separable in $G$} if,
 given any $g \in G\setminus H$, there exists a subgroup $K < G$ of finite
 index with $H < K$ and $g\notin K$. $G$ is called {\em subgroup separable}
 (or {\em LERF}) if   all finitely generated subgroups of $G$ are
 separable in  $G$.
 The profinite topology on a group $G$ is defined by proclaiming all
 finite index subgroups of $G$ to be a basis of open neighbourhoods of
 the identity.
 Since open subgroups are closed in the profinite topology, the following
 reformulates separability:

 \begin{lemma}
 \label{profinite}
 Let $G$ be a group and $H < G$ is a subgroup.
 $G$ is $H$-subgroup separable if and only if $H$ is
 closed in the profinite topology on $G$.\qed\end{lemma}

 \begin{lemma}
 \label{supergroup}
 Let $G$ be a group and $H <  K < G$. Assume that $H$ is separable in $G$ and
 that $[K:H] < \infty$. Then $K$ is separable in $G$.
 \end{lemma}
 \demo By Lemma \ref{profinite}, $H$ is closed in the profinite topology
 on $G$. Standard properties of topological groups imply that any coset
 $gH$ of $H$ in $G$ is therefore a closed subset. Since $[K:H] < \infty$,
 $K$ is a finite union of closed sets, hence closed, and therefore
 separable in $G$.\qed\\
 \vspace{2pt}
 The following is also well-known (see \cite{Lo}):
 \begin{lemma}
 \label{abelian}
 Let $G$ be a residually finite group, and $A$ a maximal abelian subgroup.
 Then $A$ is separable in $G$.\qed\end{lemma}

 \section{Proof of Theorems \ref{introgromov} and \ref{introflatinject}}
 This section is devoted to proving theorems \ref{introgromov} and \ref{introflatinject}.  We prove the
latter first, which we restate for the reader's convenience:
 
\begin{theoremrepeat} 
Let $\Gamma$ be the fundamental group  of a flat $n$-manifold.
 Then there is a quadratic form $q_{n+2}$ defined over $\bf Q$, of
 signature $(n+1,1)$ for which $\Gamma$ embeds as a subgroup of
 $\O_0(q_{n+2};{\bf Z})$. \end{theoremrepeat}

Before embarking on the proof we remark that the first part of the proof
can be replaced by the argument in the proof of Bieberbach's third
theorem giving an integral representation into $\GL(n+1, {\bf Z})$
of a Bieberbach group. However, we will use some additional features of
the construction given below in completing the proof of Theorem
\ref{introgromov}.

\proof[Proof of Theorem \ref{introflatinject}]
 Suppose that $\Gamma$ is the fundamental group of a flat $n$-manifold,
 so as discussed in \S 2.2, we have a free action of $\Gamma$
 on ${\bf E}^n$ by rigid motions.  Thus, if $g \in \Gamma$, then $g$ acts as
 $$v \rightarrow \theta(g)v + t_g.$$
 and the assignment $$g
 \rightarrow \theta(g)$$ is a homomorphism of $\Gamma$ to its holonomy group
 $\theta(\Gamma)$, with kernel the maximal translation subgroup of
 $\Gamma$.

 Suppose that $\mu_1,....,\mu_n$ generate the maximal normal free abelian
 ${\bf Z}^n$ in $G$, where $\mu_i$ acts as translation by ${\bf m}_i$,
 where we declare that this is the vector ${\bf m}_i =(0,...,1,...,0)$, one in
 $i-th$ place.  The group $\Gamma$ acts by conjugacy on the subgroup
 $\langle\mu_1,....,\mu_n\rangle$  and a calculation reveals that $g\mu_ig^{-1}$
is the
 translation given by $$v \rightarrow v + \theta(g){\bf m}_i.$$
The  normality of the translation subgroup shows that $$g\mu_ig^{-1} =
 \mu_1^{a_{i,1}(g)}......\mu_n^{a_{i,k}(g)}$$ for some collection of
 integers $\{ a_{i,j}(g) \}$.  Equating these two statements gives a
 finite integral representation of $G$ given by $$ \theta(g){\bf m}_i
 = \Sigma_j a_{i,j}(g) {\bf m}_j $$ We now construct an integral
 linear representation of $\Gamma$, as follows. Choose a presentation
 for the group $\Gamma$ using generators $g_1,.....,g_p$ and with relations
 $w_t(g_1,...,g_p) = I$, and add all relators which say
 $w_j(g_1,...,g_p) = \mu_j$ for each $1 \leq j \leq n$. (These ensure
 that the chosen ${\bf m}_i$'s don't change.) Each $g_i$ acts as $$v
 \rightarrow \theta(g_i)v + t_i,$$ so that expanding the equations
 coming from the relators, we get a collection of equations for the
 $t_i$'s with rational coefficients which have some solution (for
 example, that coming from the identity representation that we are
 given for $\Gamma$ as a Bieberbach group).

 It follows that there are rational solutions to these equations and
 we claim that any such solution gives a {\em faithful} and rational
 representation of $G$.  Pick any rational solution and regard this as
 a representation $\rho : \Gamma \rightarrow \rho(\Gamma)$. The
 conditions imposed by the second batch of equations guarantee that
 the restriction of $\rho$ to the translation subgroup of $\Gamma$ is
 actually the identity homomorphism.  Since the translation subgroup
 is isomorphic to ${\bf Z}^n$ and this is Hopfian, it follows that
 $ker(\rho)$ avoids the translation subgroup of finite index and hence
 $ker(\rho)$ is trivial, since $\Gamma$ is torsion free.  This shows
 that $\rho$ is an isomorphism as required.

 It follows that there is a {\em rational} solution for the $t_i$'s in
 terms of the ${\bf m}_i$'s. Convert the affine representation on
 ${\bf E}^n$ to a rational linear representation on ${\bf E}^{n+1}$ by
 $$\Phi^*(g) = \left( \begin{array}{ccc} \theta(g) &| & t_g \\ 0 & |
 & 1 \end{array} \right).$$
This is a faithful rational linear
 representation of $\Gamma$, which comes from coning the given action
 of $\Gamma$ by rigid motions in the hyperplane $e_{n+1} = 1$ to the
 origin in ${\bf E}^{n+1}$. By conjugating the representation, we may
 rescale the vector $e_{n+1}$ and thus arrange that the representation
 $\Phi^*$ is actually by integral matrices.

 It is slightly more convenient at this stage to define a new faithful
 integral representation by setting $$\Phi(g) = (\Phi^*(g)^T)^{-1} = \left(
 \begin{array}{ccc} \theta(g^{-1})^T &| & 0 \\ -\theta(g^{-1})(t_g)^T
 & | & 1
 \end{array} \right),$$ where $A^T$ denotes transpose. We can then extend
 this representation to ${\bf E}^{n+2}$ by mapping $g$ to  
 $$\hat{\Phi}(g) =
 \left( \begin{array}{ccc}                \Phi(g) &| & {\bf v}_g \\      
            
 0     & |  & 1              
 \end{array} \right) $$
 where the column vector ${\bf v}_g$ is to be determined.

 Now let $\langle,\rangle$ be any $\theta^T$ invariant positive definite inner
 product on the ${\bf Z}$-module $\langle\mu_1,....,\mu_n\rangle$; such an inner
 product exists by taking a random inner product and forming the
 $\theta$-average. Let $D$ be the symmetric rational matrix associated
 to this form in the $\{{\bf m}_i\}$ basis.  Extend this form to ${\bf
 E}^{n+2}$ by summing on a subspace $H_2$, which in the language of
 quadratic forms is a hyperbolic plane. More precisely, we let $H_2$
 denote the 2-dimensional form $2XY$, with associated symmetric matrix
 $\pmatrix{0 & 1\cr 1 & 0}$ (see \cite{La} Chapter 1).  The form
 $D\oplus H_2$ now has signature $(n+1,1)$.

 Denoting vectors lying in ${\bf E}^n$ by $w$ and the last two
 dimensions by $v_1$ and $v_2$, it is a simple matter of linear algebra
 to show that  ${\bf v}_g = (W_g, \tau_g) \in {\bf E}^n \oplus \langle v_1\rangle$ 
 may be chosen so that each $\hat{\Phi}(g)$ is an  isometry of the form $D\oplus
H_2$.
 
The linear algebra suggests that  the  matrices $\hat{\Phi}(g)$ may
 be nonintegral in the last column, since the initial solution vectors ${\bf
v}_g$ need
 only be rational. However, conjugating by a matrix of the form  $$ \left(
 \begin{array}{ccc} {\bf I} &| & {\bf 0} \\ 0 & | & K \end{array}
 \right) $$
 we may find a new collection of matrices which are the same save in the
 last column, and which has ${\bf v}_g$ replaced by $K.{\bf v}_g$.  In particular,
 for suitable $K$ we may arrange that the conjugated representation is integral.
 
 After this conjugation, the new
 representation now leaves invariant a different form, but this new form is
 rationally equivalent to $D\oplus H_2$; in particular, it continues to
 have signature $(n+1,1)$.\\
 \vspace{4pt}
 {\bf Claim}\qua With
 this choice, we get a faithful integral
 representation of the group $\Gamma$.
 \vspace{-2pt}
 \proof We need only show that the relations in $\Gamma$ hold. Faithfulness
 will follow, because if a product of these matrices is the identity, then it
 must at least be the identity in the $n+1$ representation which is already
 a faithful representation of $\Gamma$.

 We prove the claim by showing that any isometry, $\gamma$ say, which is
the identity on the upper left $n+1 \times n+1$ block is in fact the identity.

Note that by construction, every $\gamma \in \Gamma$ stabilizes 
 $v_{1}$.

 Pick a random $w \in {\bf E}^n$. Then $0 = \langle w, v_2\rangle =
\langle\gamma w ,
 \gamma v_2\rangle = \langle w , w' + \xi v_1 + v_2\rangle = \langle w ,
w'\rangle$. This holds
 for all
 $w$, so that $w' = 0 $.

 Also $0 = \langle v_2, v_2\rangle = \langle\gamma v_2 , \gamma v_2 \rangle =
\langle \xi v_1 +
 v_2 ,
 \xi v_1 + v_2\rangle = 2\xi$ so that $\xi = 0$, implying $\gamma v_2 = v_2$.
 as required. 

This completes the proof of theorem \ref{introflatinject}. \qed

 \vspace{2pt}
 {\bf Remark}\qua Note that the
 construction exhibits $\Gamma$ explicitly as a subgroup of the stabiliser
 of the lightlike vector $v_1$.

\proof[Proof of Theorem \ref{introgromov}]
To complete the proof of Theorem \ref{introgromov} we proceed as
follows.  Let $q_{n+2}$ be the form constructed above and consider
$\O_0(q_{n+2};{\bf Z})$.  As noted in the Remark above, the
construction yields $\Gamma$ as a subgroup of the stabiliser in
$\O_0(q_{n+2};{\bf Z})$ of
the light-like vector $v_1$, however the proof provides no control
over whether $\Gamma$ is actually equal to $\stab(v_1)$.  To achieve this,
we show that the subgroup $\Gamma$ is separable in $\O_0(q_{n+2};{\bf Z})$ and
then the
theorem follows by a standard separability argument.

To this end, let $C$ be the maximal peripheral subgroup of
$\O_0(q_{n+2};{\bf Z})$ fixing $v_1$, so that $\Gamma < C <
\O_0(q_{k+2};{\bf Z})$, and $C$ is a crystallographic group.  We
recall that by Theorem \ref{crystal}, $C$ contains a translational
subgroup $T^*$ which is free abelian of rank $n$, and is the maximal
abelian subgroup of $C$.
  
We begin by observing that our construction of the group $\Gamma$ began
with a translational subgroup which contained translation by $1$ in
each of the coordinate directions, so that after the two dilation
conjugacies which convert rational to integral, we see that for some
integer $p \rangle 1 $, the maximal translational subgroup $T^{*}$ of
$\stab(v_{1}) \leq \O_0(q_{n+2};{\bf Z})$ contains the group $T_{p}$
consisting of translations by $p$ in each of the coordinate directions
of the first $ n + 1$ coordinates.
  
We claim that $T_{p}$ is separable in $\O_0(q_{n+2};{\bf Z})$.  The
reason is this: Firstly, any element of $\O_0(q_{n+2};{\bf Z})$ which
lies outside $T^{*}$ can be separated from $T_{p}$ since in fact it
can be separated from $T^{*}$ by Lemma \ref{abelian}.
  
Secondly, we claim that any element of $T^{*} - T_{p}$ may be
separated from $T_{p}$.  This involves a few cases, which we now
sketch.
  
In the first place, we observe that all the elements of $\gamma \in T^{*}$ 
have the first $n+1$ entries of its last
row being zero, since $\langle\gamma(w), v_{1}\rangle = \langle w ,
\gamma^{-1}(v_{1})\rangle = 
\langle w , v_{1}\rangle = 0$.

Moreover, if we look at the upper left $n \times n$ block of any element
of $T_{p}$, this is constructed to be the identity matrix and if the
element $\gamma$ we wish to separate does not have this property then
we may separate by choosing a random prime $q$ not dividing some
such entry and use the restriction of the homomorphism $\SL(n+2,{\bf
Z}) \rightarrow \SL(n+2,{\bf Z}/q{\bf Z})$ to the integral subgroup
$\O_0(q_{n+2};{\bf Z})$. 

It follows that it remains to separate an element $\gamma \in T^{*} - T_{p}$
which is the identity matrix in the first $n+1$ columns save for the first
$n$ entries in the $n+1$-st row. It is these  entries which
contribute to the translational nature of the elements of $T_{p}$.
However, recall that we have proved that any isometry of $\langle,\rangle$ which
is the identity on the upper left $n+1 \times n+1$ block must in fact be
the identity matrix.  It follows that in the matrix $\gamma$ these
entries cannot all be divisible by $p$ (else the upper left $n+1 \times n+1$
block is identical with that for some matrix of $T_{p}$ and we deduce that
$\gamma \in T_{p}$) so that we may use the reduction map $\SL(n+2,{\bf
Z}) \rightarrow \SL(n+2,{\bf Z}/p{\bf Z})$ to separate $\gamma$ in
this case.
  
We note that this argument can actually be used to show a little more,
namely that for any integer $r$, the subgroups $T_{rp}$ of $T_{p}$ are
separable in $\O_0(q_{n+2};{\bf Z})$.
  
The separability of $\Gamma$ may now be deduced.  For if we let $T_{\Gamma}$ be
the maximal abelian subgroup of $\Gamma$, then $T_{\Gamma} \cap T_{p}$ is a
subgroup of finite index in $T_{p}$ and it follows that there is an
integer $r$ for which $T_{rp} \leq T_{\Gamma} \cap T_{p} \leq \Gamma$.  The
separability of $\Gamma$ follows from Lemma \ref{supergroup}.
This completes the proof of theorem \ref{introgromov}. \qed

\vspace{2pt}
 {\bf Remark}\qua
 The number of cusps for the groups $\O_0(q_{n+2};{\bf Z})$ can be greater
 than one, even for simple examples. For example, the groups
 $\O_0(f_n;{\bf Z})$ have 1 cusp for $2\leq n\leq 8$, but $\O_0(f_9;{\bf Z})$
 has 2 cusps.  This can be seen from \cite{Vi} and \cite{VS} which
 describes these unit groups as groups generated by reflections
 in certain ideal simplices in ${\bf H}^n$. The number of cusps
 being easy to read off from the Coxeter diagrams.

 \section{Example}

 We finish off by giving an example of the construction as an aid to
 the proof of Theorem \ref{introflatinject}.  We shall take the
 Hantsche-Wendt manifold, which arises as the 3-fold cyclic branched
 cover of the figure-eight knot. Its fundamental group
 is therefore a Bieberbach group in dimension
 $3$. Representing matrices are provided on pp. 6-7 of \cite{Ch} which
 in our notation are:
 $$\Phi^*(a) =
  \pmatrix{ -1 & 0 &0 &1/2 \cr
             0 & -1& 0 & 1/2 \cr
             0& 0 &1 &1/2\cr
             0&0&0&1}
 \;\;\;\hbox{and}\;\;\;
 \Phi^*(b) =  \pmatrix{ 1 & 0 &0 &1/2 \cr
             0 & -1& 0 & 0 \cr
             0& 0 &-1 &0\cr
             0&0&0&1}$$
 We conjugate to clear fractions and form the representation $\Phi$ and hence
 $\hat{\Phi}$.
 In this case an invariant form for the finite holonomy is $D=x^2+y^2+z^2$ and
 when we
 solve the equations for the group to be an isometry for the form $D\oplus H_2$
 we obtain:
 $$\hat{\Phi} (a) = \pmatrix{ -1&0&0&0&1\cr
                              0&-1&0&0&1\cr
                                0&0&1&0&1\cr
                               1&1&-1&1&-{3\over2} \cr
                               0&0&0&0&1}\;\;\;\hbox{and}\;\;\;
 \hat{\Phi} (b) =
 \pmatrix{ 1&0&0&0&1\cr
            0&-1&0&0&0\cr
             0&0&-1&0&0\cr
            -1&0&0&1&-{1\over2} \cr
             0&0&0&0&1}$$
 Letting $K=2$, and conjugating, gives integral matrices,
 $$\hat{\Phi} (a) = \pmatrix{ -1&0&0&0&2\cr
                              0&-1&0&0&2\cr
                                0&0&1&0&2\cr
                               1&1&-1&1&-3 \cr
                               0&0&0&0&1}\;\;\;\hbox{and}\;\;\;
 \hat{\Phi} (b) =
 \pmatrix{ 1&0&0&0&2\cr
            0&-1&0&0&0\cr
             0&0&-1&0&0\cr
            -1&0&0&1&-1 \cr
             0&0&0&0&1}$$
 preserving the rationally equivalent form $x^2 + y^2 + z^2 + 4wt$.

 Note that the above form is equivalent over $\bf Q$ to
 $f_4 = x^2 + y^2 + z^2 + w^2 - t^2$. In \cite{RT} the authors obtain
 the Hantsche-Wendt manifold as a cusp cross-section of a hyperbolic 4-manifold
 arising from a torsion-free subgroup in $O_0(f_4;{\bf Z})$.

\vspace{4pt}
{\bf Acknowledgements}\qua The first author was partially supported
 by the N.S.F, and the second was partially supported by the N.S.F, the
 Alfred P. Sloan Foundation and a grant from the Texas Advanced
 Research Program.

\Addresses\recd
\end{document}

%% file: agtout.tex
\def\ifplaintex{\expandafter\ifx\csname documentclass\endcsname\relax}

\def\gtp{{\mathsurround=0pt\it $\cal G\mskip-2mu$eometry \&\ 
$\cal T\!\!$opology $\cal P\!$ublications}}  

\def\recd{{\small Received:\qua\receiveddate\ifx\reviseddate\relax
\else\qquad Revised:\qua\reviseddate\fi\par}} 

\def\lognumber#1{\def\thelognumber{#1}}
\def\volumenumber#1{\def\thevolumenumber{#1}}
\def\volumeyear#1{\def\thevolumeyear{#1}}
\def\papernumber#1{\def\thepapernumber{#1}}
\def\pagenumbers#1#2{\def\startpage{#1}\def\finishpage{#2}}
\def\published#1{\def\publishdate{#1}}

\def\received#1{\def\receiveddate{#1}}

\def\accepted#1{\def\accepteddate{#1}}

\def\asciiaddress#1{\def\theasciiaddress{#1}}

\long\def\asciiabstract#1{\long\def\theasciiabstract{#1}}

\let\\\par\let\thelognumber\relax\let\thevolumenumber\relax
\let\thepapernumber\relax\let\thevolumeyear\relax\let\startpage\relax
\let\finishpage\relax\let\publishdate\relax\let\receiveddate\relax
\let\reviseddate\relax\let\accepteddate\relax\let\theasciititle\relax
\let\theasciiauthors\relax\let\theasciiaddress\relax
\let\theasciiabstract\relax

\let\theasciiemail\relax

\ifplaintex
\font\logobig=cmssbx10 scaled 3836
\font\logomed=cmssbx10 scaled 2557
\else
\font\logobig=cmssbx10 scaled 4200
\font\logomed=cmssbx10 scaled 2800
\fi

\long\def\makeagttitle{   
\count0=\startpage
\agt\hfill      
\hbox to 45truept{\vbox to 0pt{\vglue -13truept{\logomed A\kern -.37em{\logobig 
T}\kern -.38em G}\vss}\hss}
\break
{\small Volume \thevolumenumber\ (\thevolumeyear)
\startpage--\finishpage\nl
Published: \publishdate}

\vglue .25truein

{\parskip=0pt\leftskip 0pt plus
1fil\def\\{\par\smallskip}{\Large\bf\thetitle}\par\medskip} \vglue
0.05truein

{\parskip=0pt\leftskip 0pt plus 1fil\def\\{\par}{\sc\theauthors}
\par\medskip}%
 
\vglue 0.03truein 

{\small\leftskip 25truept\rightskip 25truept{\bf Abstract}\stdspace\theabstract

{\bf AMS Classification}\stdspace\theprimaryclass
\ifx\thesecondaryclass\relax\else; \thesecondaryclass\fi\par
{\bf Keywords}\stdspace \thekeywords\par}\vglue 7truept

}   

\ifplaintex
\font\phead=cmsl9 scaled 950
\font\pnum=cmbx10 scaled 913
\font\pfoot=cmsl9 scaled 950
\headline{\vbox to 0pt{\vskip -4.5mm\line{\small\phead\ifnum
\count0=\startpage ISSN 1472-2739 (on-line) 1472-2747 (printed)
\hfill {\pnum\folio}\else\ifodd\count0\def\\{ }%
\ifx\theshorttitle\relax\thetitle\else\theshorttitle\fi\hfill{\pnum\folio}
\else\def\\{ and }{\pnum\folio}\hfill\ifx\theshortauthors\relax\theauthors
\else\theshortauthors\fi\fi\fi}\vss}}
\footline{\vbox to 0pt{\vglue 0mm\line{\small\pfoot\ifnum\count0=\startpage
\copyright\ \gtp\hfill\else
\agt, Volume \thevolumenumber\ (\thevolumeyear)\hfill\fi}\vss}}
\else
\font=cmsl9 scaled 1050
\font\lnum=cmbx10 
\font=cmsl9 scaled 1050
\makeatletter
\def\@oddhead{{\small\lhead\ifnum\count0=\startpage ISSN 1472-2739 
(on-line) 1472-2747 (printed)\hfill {\lnum\number\count0}\else\ifodd\count0
\def\\{ }\ifx\theshorttitle\relax \thetitle \else\theshorttitle\fi\hfill
{\lnum\number\count0}\else\def\\{ and }{\lnum\number\count0}
\hfill\ifx\theshortauthors\relax 
\theauthors\else\theshortauthors\fi\fi\fi}}\def\@evenhead{\@oddhead}
\def\@oddfoot{\small\lfoot\ifnum\count0=\startpage\copyright\ \gtp\hfill\else
\agt, Volume \thevolumenumber\ (\thevolumeyear)\hfill\fi}
\def\@evenfoot{\@oddfoot}
\makeatother
\fi
\let\maketitlepage\makeagttitle
\let\makeshorttitle\maketitlepage
\let\maketitle\maketitlepage

\newwrite\gtoutfile
\long\gdef\makeheadfile{  
{\def\\{, }\def\s{ }
\immediate\openout\gtoutfile head.xxx
\immediate\write\gtoutfile{To: math@arxiv.org}
\immediate\write\gtoutfile{Subject: put OR rep NNNNN:ppppp}
\immediate\write\gtoutfile{--text follows this line--}
\immediate\write\gtoutfile{Proxy-for: \ifx\theasciiauthors\relax
\theauthors\else\theasciiauthors\fi\s<\ifx\theasciiemail\relax\theemail\else\theasciiemail\fi>}
\immediate\write\gtoutfile{\noexpand\\}
\immediate\write\gtoutfile{Authors: \ifx\theasciiauthors\relax
\theauthors\else\theasciiauthors\fi}
{\def\\{ }\immediate\write\gtoutfile{Title: \ifx\theasciititle\relax
\thetitle\else\theasciititle\fi}}
\immediate\write\gtoutfile{Subj-class: GT or SG, GR etc}
\immediate\write\gtoutfile{MSC-class: \theprimaryclass\ifx\thesecondaryclass\relax\else, \thesecondaryclass\fi}
\immediate\write\gtoutfile{Journal-ref: Algebr. Geom. Topol. \thevolumenumber\s
(\thevolumeyear) \startpage-\finishpage}
\immediate\write\gtoutfile{Comments: Published by Algebraic and
Geometric Topology at}
\immediate\write\gtoutfile{\s\s\s  http://www.maths.warwick.ac.uk/agt/AGTVol\thevolumenumber/agt-\thevolumenumber-\thepapernumber.abs.html}
\immediate\write\gtoutfile{\noexpand\\}
\immediate\write\gtoutfile{}
\ifx\theasciiabstract\relax
\immediate\write\gtoutfile{\theabstract}\else
\immediate\write\gtoutfile{\theasciiabstract}\fi
\immediate\write\gtoutfile{}
\immediate\write\gtoutfile{\noexpand\\}
\immediate\write\gtoutfile{}
\immediate\closeout\gtoutfile}}  

\def\maketitlepage{\makeagttitle\makeheadfile}
\let\makeshorttitle\maketitlepage
\let\maketitle\maketitlepage

\def\ifplaintex{\expandafter\ifx\csname documentclass\endcsname\relax}

\def\gtp{{\mathsurround=0pt\it $\cal G\mskip-2mu$eometry \&\ 
$\cal T\!\!$opology $\cal P\!$ublications}}  

\def\recd{{\small Received:\qua\receiveddate\ifx\reviseddate\relax
\else\qquad Revised:\qua\reviseddate\fi\par}} 

\def\lognumber#1{\def\thelognumber{#1}}
\def\volumenumber#1{\def\thevolumenumber{#1}}
\def\volumeyear#1{\def\thevolumeyear{#1}}
\def\papernumber#1{\def\thepapernumber{#1}}
\def\pagenumbers#1#2{\def\startpage{#1}\def\finishpage{#2}}
\def\published#1{\def\publishdate{#1}}

\def\received#1{\def\receiveddate{#1}}

\def\accepted#1{\def\accepteddate{#1}}

\def\asciiaddress#1{\def\theasciiaddress{#1}}

\long\def\asciiabstract#1{\long\def\theasciiabstract{#1}}

\let\\\par\let\thelognumber\relax\let\thevolumenumber\relax
\let\thepapernumber\relax\let\thevolumeyear\relax\let\startpage\relax
\let\finishpage\relax\let\publishdate\relax\let\receiveddate\relax
\let\reviseddate\relax\let\accepteddate\relax\let\theasciititle\relax
\let\theasciiauthors\relax\let\theasciiaddress\relax
\let\theasciiabstract\relax

\let\theasciiemail\relax

\ifplaintex
\font\logobig=cmssbx10 scaled 3836
\font\logomed=cmssbx10 scaled 2557
\else
\font\logobig=cmssbx10 scaled 4200
\font\logomed=cmssbx10 scaled 2800
\fi

\long\def\makeagttitle{   
\count0=\startpage
\agt\hfill      
\hbox to 45truept{\vbox to 0pt{\vglue -13truept{\logomed A\kern -.37em{\logobig 
T}\kern -.38em G}\vss}\hss}
\break
{\small Volume \thevolumenumber\ (\thevolumeyear)
\startpage--\finishpage\nl
Published: \publishdate}

\vglue .25truein

{\parskip=0pt\leftskip 0pt plus
1fil\def\\{\par\smallskip}{\Large\bf\thetitle}\par\medskip} \vglue
0.05truein

{\parskip=0pt\leftskip 0pt plus 1fil\def\\{\par}{\sc\theauthors}
\par\medskip}%
 
\vglue 0.03truein 

{\small\leftskip 25truept\rightskip 25truept{\bf Abstract}\stdspace\theabstract

{\bf AMS Classification}\stdspace\theprimaryclass
\ifx\thesecondaryclass\relax\else; \thesecondaryclass\fi\par
{\bf Keywords}\stdspace \thekeywords\par}\vglue 7truept

}   

\ifplaintex
\font\phead=cmsl9 scaled 950
\font\pnum=cmbx10 scaled 913
\font\pfoot=cmsl9 scaled 950
\headline{\vbox to 0pt{\vskip -4.5mm\line{\small\phead\ifnum
\count0=\startpage ISSN 1472-2739 (on-line) 1472-2747 (printed)
\hfill {\pnum\folio}\else\ifodd\count0\def\\{ }%
\ifx\theshorttitle\relax\thetitle\else\theshorttitle\fi\hfill{\pnum\folio}
\else\def\\{ and }{\pnum\folio}\hfill\ifx\theshortauthors\relax\theauthors
\else\theshortauthors\fi\fi\fi}\vss}}
\footline{\vbox to 0pt{\vglue 0mm\line{\small\pfoot\ifnum\count0=\startpage
\copyright\ \gtp\hfill\else
\agt, Volume \thevolumenumber\ (\thevolumeyear)\hfill\fi}\vss}}
\else
\font=cmsl9 scaled 1050
\font\lnum=cmbx10 
\font=cmsl9 scaled 1050
\makeatletter
\def\@oddhead{{\small\lhead\ifnum\count0=\startpage ISSN 1472-2739 
(on-line) 1472-2747 (printed)\hfill {\lnum\number\count0}\else\ifodd\count0
\def\\{ }\ifx\theshorttitle\relax \thetitle \else\theshorttitle\fi\hfill
{\lnum\number\count0}\else\def\\{ and }{\lnum\number\count0}
\hfill\ifx\theshortauthors\relax 
\theauthors\else\theshortauthors\fi\fi\fi}}\def\@evenhead{\@oddhead}
\def\@oddfoot{\small\lfoot\ifnum\count0=\startpage\copyright\ \gtp\hfill\else
\agt, Volume \thevolumenumber\ (\thevolumeyear)\hfill\fi}
\def\@evenfoot{\@oddfoot}
\makeatother
\fi
\let\maketitlepage\makeagttitle
\let\makeshorttitle\maketitlepage
\let\maketitle\maketitlepage

\newwrite\gtoutfile
\long\gdef\makeheadfile{  
{\def\\{, }\def\s{ }
\immediate\openout\gtoutfile head.xxx
\immediate\write\gtoutfile{To: math@arxiv.org}
\immediate\write\gtoutfile{Subject: put OR rep NNNNN:ppppp}
\immediate\write\gtoutfile{--text follows this line--}
\immediate\write\gtoutfile{Proxy-for: \ifx\theasciiauthors\relax
\theauthors\else\theasciiauthors\fi\s<\ifx\theasciiemail\relax\theemail\else\theasciiemail\fi>}
\immediate\write\gtoutfile{\noexpand\\}
\immediate\write\gtoutfile{Authors: \ifx\theasciiauthors\relax
\theauthors\else\theasciiauthors\fi}
{\def\\{ }\immediate\write\gtoutfile{Title: \ifx\theasciititle\relax
\thetitle\else\theasciititle\fi}}
\immediate\write\gtoutfile{Subj-class: GT or SG, GR etc}
\immediate\write\gtoutfile{MSC-class: \theprimaryclass\ifx\thesecondaryclass\relax\else, \thesecondaryclass\fi}
\immediate\write\gtoutfile{Journal-ref: Algebr. Geom. Topol. \thevolumenumber\s
(\thevolumeyear) \startpage-\finishpage}
\immediate\write\gtoutfile{Comments: Published by Algebraic and
Geometric Topology at}
\immediate\write\gtoutfile{\s\s\s  http://www.maths.warwick.ac.uk/agt/AGTVol\thevolumenumber/agt-\thevolumenumber-\thepapernumber.abs.html}
\immediate\write\gtoutfile{\noexpand\\}
\immediate\write\gtoutfile{}
\ifx\theasciiabstract\relax
\immediate\write\gtoutfile{\theabstract}\else
\immediate\write\gtoutfile{\theasciiabstract}\fi
\immediate\write\gtoutfile{}
\immediate\write\gtoutfile{\noexpand\\}
\immediate\write\gtoutfile{}
\immediate\closeout\gtoutfile}}  

\def\maketitlepage{\makeagttitle\makeheadfile}
\let\makeshorttitle\maketitlepage
\let\maketitle\maketitlepage

\def\ifplaintex{\expandafter\ifx\csname documentclass\endcsname\relax}

\def\gtp{{\mathsurround=0pt\it $\cal G\mskip-2mu$eometry \&\ 
$\cal T\!\!$opology $\cal P\!$ublications}}  

\def\recd{{\small Received:\qua\receiveddate\ifx\reviseddate\relax
\else\qquad Revised:\qua\reviseddate\fi\par}} 

\def\lognumber#1{\def\thelognumber{#1}}
\def\volumenumber#1{\def\thevolumenumber{#1}}
\def\volumeyear#1{\def\thevolumeyear{#1}}
\def\papernumber#1{\def\thepapernumber{#1}}
\def\pagenumbers#1#2{\def\startpage{#1}\def\finishpage{#2}}
\def\published#1{\def\publishdate{#1}}

\def\received#1{\def\receiveddate{#1}}

\def\accepted#1{\def\accepteddate{#1}}

\def\asciiaddress#1{\def\theasciiaddress{#1}}

\long\def\asciiabstract#1{\long\def\theasciiabstract{#1}}

\let\\\par\let\thelognumber\relax\let\thevolumenumber\relax
\let\thepapernumber\relax\let\thevolumeyear\relax\let\startpage\relax
\let\finishpage\relax\let\publishdate\relax\let\receiveddate\relax
\let\reviseddate\relax\let\accepteddate\relax\let\theasciititle\relax
\let\theasciiauthors\relax\let\theasciiaddress\relax
\let\theasciiabstract\relax

\let\theasciiemail\relax

\ifplaintex
\font\logobig=cmssbx10 scaled 3836
\font\logomed=cmssbx10 scaled 2557
\else
\font\logobig=cmssbx10 scaled 4200
\font\logomed=cmssbx10 scaled 2800
\fi

\long\def\makeagttitle{   
\count0=\startpage
\agt\hfill      
\hbox to 45truept{\vbox to 0pt{\vglue -13truept{\logomed A\kern -.37em{\logobig 
T}\kern -.38em G}\vss}\hss}
\break
{\small Volume \thevolumenumber\ (\thevolumeyear)
\startpage--\finishpage\nl
Published: \publishdate}

\vglue .25truein

{\parskip=0pt\leftskip 0pt plus
1fil\def\\{\par\smallskip}{\Large\bf\thetitle}\par\medskip} \vglue
0.05truein

{\parskip=0pt\leftskip 0pt plus 1fil\def\\{\par}{\sc\theauthors}
\par\medskip}%
 
\vglue 0.03truein 

{\small\leftskip 25truept\rightskip 25truept{\bf Abstract}\stdspace\theabstract

{\bf AMS Classification}\stdspace\theprimaryclass
\ifx\thesecondaryclass\relax\else; \thesecondaryclass\fi\par
{\bf Keywords}\stdspace \thekeywords\par}\vglue 7truept

}   

\ifplaintex
\font\phead=cmsl9 scaled 950
\font\pnum=cmbx10 scaled 913
\font\pfoot=cmsl9 scaled 950
\headline{\vbox to 0pt{\vskip -4.5mm\line{\small\phead\ifnum
\count0=\startpage ISSN 1472-2739 (on-line) 1472-2747 (printed)
\hfill {\pnum\folio}\else\ifodd\count0\def\\{ }%
\ifx\theshorttitle\relax\thetitle\else\theshorttitle\fi\hfill{\pnum\folio}
\else\def\\{ and }{\pnum\folio}\hfill\ifx\theshortauthors\relax\theauthors
\else\theshortauthors\fi\fi\fi}\vss}}
\footline{\vbox to 0pt{\vglue 0mm\line{\small\pfoot\ifnum\count0=\startpage
\copyright\ \gtp\hfill\else
\agt, Volume \thevolumenumber\ (\thevolumeyear)\hfill\fi}\vss}}
\else
\font=cmsl9 scaled 1050
\font\lnum=cmbx10 
\font=cmsl9 scaled 1050
\makeatletter
\def\@oddhead{{\small\lhead\ifnum\count0=\startpage ISSN 1472-2739 
(on-line) 1472-2747 (printed)\hfill {\lnum\number\count0}\else\ifodd\count0
\def\\{ }\ifx\theshorttitle\relax \thetitle \else\theshorttitle\fi\hfill
{\lnum\number\count0}\else\def\\{ and }{\lnum\number\count0}
\hfill\ifx\theshortauthors\relax 
\theauthors\else\theshortauthors\fi\fi\fi}}\def\@evenhead{\@oddhead}
\def\@oddfoot{\small\lfoot\ifnum\count0=\startpage\copyright\ \gtp\hfill\else
\agt, Volume \thevolumenumber\ (\thevolumeyear)\hfill\fi}
\def\@evenfoot{\@oddfoot}
\makeatother
\fi
\let\maketitlepage\makeagttitle
\let\makeshorttitle\maketitlepage
\let\maketitle\maketitlepage

\newwrite\gtoutfile
\long\gdef\makeheadfile{  
{\def\\{, }\def\s{ }
\immediate\openout\gtoutfile head.xxx
\immediate\write\gtoutfile{To: math@arxiv.org}
\immediate\write\gtoutfile{Subject: put OR rep NNNNN:ppppp}
\immediate\write\gtoutfile{--text follows this line--}
\immediate\write\gtoutfile{Proxy-for: \ifx\theasciiauthors\relax
\theauthors\else\theasciiauthors\fi\s<\ifx\theasciiemail\relax\theemail\else\theasciiemail\fi>}
\immediate\write\gtoutfile{\noexpand\\}
\immediate\write\gtoutfile{Authors: \ifx\theasciiauthors\relax
\theauthors\else\theasciiauthors\fi}
{\def\\{ }\immediate\write\gtoutfile{Title: \ifx\theasciititle\relax
\thetitle\else\theasciititle\fi}}
\immediate\write\gtoutfile{Subj-class: GT or SG, GR etc}
\immediate\write\gtoutfile{MSC-class: \theprimaryclass\ifx\thesecondaryclass\relax\else, \thesecondaryclass\fi}
\immediate\write\gtoutfile{Journal-ref: Algebr. Geom. Topol. \thevolumenumber\s
(\thevolumeyear) \startpage-\finishpage}
\immediate\write\gtoutfile{Comments: Published by Algebraic and
Geometric Topology at}
\immediate\write\gtoutfile{\s\s\s  http://www.maths.warwick.ac.uk/agt/AGTVol\thevolumenumber/agt-\thevolumenumber-\thepapernumber.abs.html}
\immediate\write\gtoutfile{\noexpand\\}
\immediate\write\gtoutfile{}
\ifx\theasciiabstract\relax
\immediate\write\gtoutfile{\theabstract}\else
\immediate\write\gtoutfile{\theasciiabstract}\fi
\immediate\write\gtoutfile{}
\immediate\write\gtoutfile{\noexpand\\}
\immediate\write\gtoutfile{}
\immediate\closeout\gtoutfile}}  

\def\maketitlepage{\makeagttitle\makeheadfile}
\let\makeshorttitle\maketitlepage
\let\maketitle\maketitlepage

\def\ifplaintex{\expandafter\ifx\csname documentclass\endcsname\relax}

\def\gtp{{\mathsurround=0pt\it $\cal G\mskip-2mu$eometry \&\ 
$\cal T\!\!$opology $\cal P\!$ublications}}  

\def\recd{{\small Received:\qua\receiveddate\ifx\reviseddate\relax
\else\qquad Revised:\qua\reviseddate\fi\par}} 

\def\lognumber#1{\def\thelognumber{#1}}
\def\volumenumber#1{\def\thevolumenumber{#1}}
\def\volumeyear#1{\def\thevolumeyear{#1}}
\def\papernumber#1{\def\thepapernumber{#1}}
\def\pagenumbers#1#2{\def\startpage{#1}\def\finishpage{#2}}
\def\published#1{\def\publishdate{#1}}

\def\received#1{\def\receiveddate{#1}}

\def\accepted#1{\def\accepteddate{#1}}

\def\asciiaddress#1{\def\theasciiaddress{#1}}

\long\def\asciiabstract#1{\long\def\theasciiabstract{#1}}

\let\\\par\let\thelognumber\relax\let\thevolumenumber\relax
\let\thepapernumber\relax\let\thevolumeyear\relax\let\startpage\relax
\let\finishpage\relax\let\publishdate\relax\let\receiveddate\relax
\let\reviseddate\relax\let\accepteddate\relax\let\theasciititle\relax
\let\theasciiauthors\relax\let\theasciiaddress\relax
\let\theasciiabstract\relax

\let\theasciiemail\relax

\ifplaintex
\font\logobig=cmssbx10 scaled 3836
\font\logomed=cmssbx10 scaled 2557
\else
\font\logobig=cmssbx10 scaled 4200
\font\logomed=cmssbx10 scaled 2800
\fi

\long\def\makeagttitle{   
\count0=\startpage
\agt\hfill      
\hbox to 45truept{\vbox to 0pt{\vglue -13truept{\logomed A\kern -.37em{\logobig 
T}\kern -.38em G}\vss}\hss}
\break
{\small Volume \thevolumenumber\ (\thevolumeyear)
\startpage--\finishpage\nl
Published: \publishdate}

\vglue .25truein

{\parskip=0pt\leftskip 0pt plus
1fil\def\\{\par\smallskip}{\Large\bf\thetitle}\par\medskip} \vglue
0.05truein

{\parskip=0pt\leftskip 0pt plus 1fil\def\\{\par}{\sc\theauthors}
\par\medskip}%
 
\vglue 0.03truein 

{\small\leftskip 25truept\rightskip 25truept{\bf Abstract}\stdspace\theabstract

{\bf AMS Classification}\stdspace\theprimaryclass
\ifx\thesecondaryclass\relax\else; \thesecondaryclass\fi\par
{\bf Keywords}\stdspace \thekeywords\par}\vglue 7truept

}   

\ifplaintex
\font\phead=cmsl9 scaled 950
\font\pnum=cmbx10 scaled 913
\font\pfoot=cmsl9 scaled 950
\headline{\vbox to 0pt{\vskip -4.5mm\line{\small\phead\ifnum
\count0=\startpage ISSN 1472-2739 (on-line) 1472-2747 (printed)
\hfill {\pnum\folio}\else\ifodd\count0\def\\{ }%
\ifx\theshorttitle\relax\thetitle\else\theshorttitle\fi\hfill{\pnum\folio}
\else\def\\{ and }{\pnum\folio}\hfill\ifx\theshortauthors\relax\theauthors
\else\theshortauthors\fi\fi\fi}\vss}}
\footline{\vbox to 0pt{\vglue 0mm\line{\small\pfoot\ifnum\count0=\startpage
\copyright\ \gtp\hfill\else
\agt, Volume \thevolumenumber\ (\thevolumeyear)\hfill\fi}\vss}}
\else
\font=cmsl9 scaled 1050
\font\lnum=cmbx10 
\font=cmsl9 scaled 1050
\makeatletter
\def\@oddhead{{\small\lhead\ifnum\count0=\startpage ISSN 1472-2739 
(on-line) 1472-2747 (printed)\hfill {\lnum\number\count0}\else\ifodd\count0
\def\\{ }\ifx\theshorttitle\relax \thetitle \else\theshorttitle\fi\hfill
{\lnum\number\count0}\else\def\\{ and }{\lnum\number\count0}
\hfill\ifx\theshortauthors\relax 
\theauthors\else\theshortauthors\fi\fi\fi}}\def\@evenhead{\@oddhead}
\def\@oddfoot{\small\lfoot\ifnum\count0=\startpage\copyright\ \gtp\hfill\else
\agt, Volume \thevolumenumber\ (\thevolumeyear)\hfill\fi}
\def\@evenfoot{\@oddfoot}
\makeatother
\fi
\let\maketitlepage\makeagttitle
\let\makeshorttitle\maketitlepage
\let\maketitle\maketitlepage

\newwrite\gtoutfile
\long\gdef\makeheadfile{  
{\def\\{, }\def\s{ }
\immediate\openout\gtoutfile head.xxx
\immediate\write\gtoutfile{To: math@arxiv.org}
\immediate\write\gtoutfile{Subject: put OR rep NNNNN:ppppp}
\immediate\write\gtoutfile{--text follows this line--}
\immediate\write\gtoutfile{Proxy-for: \ifx\theasciiauthors\relax
\theauthors\else\theasciiauthors\fi\s<\ifx\theasciiemail\relax\theemail\else\theasciiemail\fi>}
\immediate\write\gtoutfile{\noexpand\\}
\immediate\write\gtoutfile{Authors: \ifx\theasciiauthors\relax
\theauthors\else\theasciiauthors\fi}
{\def\\{ }\immediate\write\gtoutfile{Title: \ifx\theasciititle\relax
\thetitle\else\theasciititle\fi}}
\immediate\write\gtoutfile{Subj-class: GT or SG, GR etc}
\immediate\write\gtoutfile{MSC-class: \theprimaryclass\ifx\thesecondaryclass\relax\else, \thesecondaryclass\fi}
\immediate\write\gtoutfile{Journal-ref: Algebr. Geom. Topol. \thevolumenumber\s
(\thevolumeyear) \startpage-\finishpage}
\immediate\write\gtoutfile{Comments: Published by Algebraic and
Geometric Topology at}
\immediate\write\gtoutfile{\s\s\s  http://www.maths.warwick.ac.uk/agt/AGTVol\thevolumenumber/agt-\thevolumenumber-\thepapernumber.abs.html}
\immediate\write\gtoutfile{\noexpand\\}
\immediate\write\gtoutfile{}
\ifx\theasciiabstract\relax
\immediate\write\gtoutfile{\theabstract}\else
\immediate\write\gtoutfile{\theasciiabstract}\fi
\immediate\write\gtoutfile{}
\immediate\write\gtoutfile{\noexpand\\}
\immediate\write\gtoutfile{}
\immediate\closeout\gtoutfile}}  

\def\maketitlepage{\makeagttitle\makeheadfile}
\let\makeshorttitle\maketitlepage
\let\maketitle\maketitlepage